\documentclass{article}

\newtheorem{theorem}{Theorem}[section]
\newtheorem{corollary}[theorem]{Corollary}

\newtheorem{lemma}[theorem]{Lemma}
\newtheorem{proposition}[theorem]{Proposition}

\evensidemargin\oddsidemargin
\begin{document}

\author{Vadim E. Levit and Eugen Mandrescu \\
Department of Computer Science\\
Holon Academic Institute of Technology\\
52 Golomb Str., P.O. Box 305\\
Holon 58102, ISRAEL\\
\{levitv, eugen\_m\}@barley.cteh.ac.il}
\title{On $\alpha ^{+}$-Stable K\"{o}nig-Egervary Graphs}
\date{}
\maketitle

\begin{abstract}
The stability number of a graph $G$, denoted by $\alpha (G)$, is the
cardinality of a stable set of maximum size in $G$. If its stability number
remains the same upon the addition of any edge, then $G$ is called $\alpha
^{+}$-stable. $G$ is a K\"{o}nig-Egervary graph if its order equals $\alpha
(G)+\mu (G)$, where $\mu (G)$ is the cardinality of a maximum matching in $G$%
. In this paper we characterize $\alpha ^{+}$-stable K\"{o}nig-Egervary
graphs, generalizing some previously known results on bipartite graphs and
trees. Namely, we prove that a K\"{o}nig-Egervary graph $G=(V,E)$ is $\alpha
^{+}$-stable if and only if either $\left| \cap \{V-S:S\in \Omega
(G)\}\right| =0$, or $\left| \cap \{V-S:S\in \Omega (G)\}\right| =1$, and $G$
has a perfect matching (where $\Omega (G)$ denotes the family of all maximum
stable sets of $G$). Using this characterization we obtain several new
findings on general K\"{o}nig-Egervary graphs, for example, the equality $%
\left| \cap \{S:S\in \Omega (G)\}\right| =\left| \cap \{V-S:S\in \Omega
(G)\}\right| $ is a necessary and sufficient condition for a
K\"{o}nig-Egervary graph $G$ to have a perfect matching.
\end{abstract}

\section{Introduction}

Throughout this paper $G=(V,E)$ is a simple (i.e., a finite, undirected,
loopless and without multiple edges) graph with vertex set $V=V(G)$ and edge
set $E=E(G).$ If $X\subset V$, then $G[X]$ is the subgraph of $G$ spanned by 
$X$. By $G-W$ we mean the subgraph $G[V-W]$, if $W\subset V(G)$. By $G-F$ we
denote the partial subgraph of $G$ obtained by deleting the edges of $F$,
for $F\subset E(G)$, and we use $G-e$, if $W$ $=\{e\}$. If $A,B$ $\subset V$
and $A\cap B=\emptyset $, then $(A,B)$ stands for the set $\{e=ab:a\in
A,b\in B,e\in E\}$.

A \textit{stable set} in $G$ is a set $A\subseteq V$ of pairwise
non-adjacent vertices. A stable set of maximum size will be referred as to a 
\textit{maximum stable set} of $G$ and its cardinality $\alpha (G)$ is the 
\textit{stability number }of $G$. Let $\Omega (G)$ stand for the set $\{S:S$ 
\textit{is a maximum stable set of} $G\}$.

A matching (i.e., a set of non-incident edges of $G$) of maximum cardinality 
$\mu (G)$ is a \textit{maximum matching}, and a \textit{perfect matching} is
one covering all the vertices of $G$. If $\left| V(G)\right| -2\left|
M\right| =1$, then $M$ is called \textit{near-perfect}, \cite{lovpl}. By $%
C_{n},K_{n},P_{n}$ we denote the chordless cycle on $n\geq $ $4$ vertices,
the complete graph on $n\geq 1$ vertices, and respectively the chordless
path on $n\geq 3$ vertices.

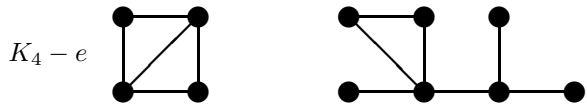
\begin{figure}[h]
\setlength{\unitlength}{1cm}%
\begin{picture}(5,1.5)\thicklines
  \multiput(4,0)(1,0){2}{\circle*{0.29}}
  \put(5,1){\circle*{0.29}}
  \put(4,1){\circle*{0.29}}
  \put(4,0){\line(1,0){1}}
  \put(4,1){\line(1,0){1}}
  \put(5,0){\line(0,1){1}} 
  \put(4,0){\line(1,1){1}}
  \put(4,0){\line(0,1){1}}
\put(3,0.5){\makebox(0,0){$K_{4}-e$}}

 \multiput(7,0)(1,0){4}{\circle*{0.29}}
 \multiput(7,1)(1,0){3}{\circle*{0.29}}
  \put(7,0){\line(1,0){3}}
  \put(7,1){\line(1,0){1}}
  \put(7,1){\line(1,-1){1}}
  \put(8,0){\line(0,1){1}} 
  \put(9,0){\line(0,1){1}} 
  
 \end{picture}
\caption{K\"{o}nig-Egervary non-bipartite graphs.}
\label{fig1}
\end{figure}

It is known that $\lfloor n/2\rfloor +1\leq \alpha (G)+\mu (G)\leq n$ holds
for any graph $G$ with $n$ vertices. Any complete graph $K_{n}$ represents
the lower bound in this inequality, while the upper bound is achieved,
according to a well-known result of Koenig, \cite{koen}, and Egervary, \cite
{eger}, by any bipartite graph. It is easy to see that there are also
non-bipartite graphs having the same property, for instance, the graphs in
Figure \ref{fig1}. 

If $\alpha (G)+\mu (G)=\left| V(G)\right| $, then $G$ is called \textit{a
K\"{o}nig-Egervary graph}. We attribute this definition to Deming \cite{dem}%
, and Sterboul \cite{ster}, but it is also possible to say that Klee \cite
{klee} defined this notion implicitly before them. These graphs were studied
by Korach \cite{kor}, Lovasz \cite{lov}, Lovasz and Plummer \cite{lovpl},
Bourjolly and Pulleyblank \cite{bourpull}, Pulleyblank \cite{pulleybl}, and
generalized by Bourjolly, Hammer and Simeone \cite{bourhams1}, Paschos and
Demange \cite{pasdema}. Since $G$ is a K\"{o}nig-Egervary graph if and only
if all its connected components are K\"{o}nig-Egervary graphs, throughout
this paper we shall consider only connected K\"{o}nig-Egervary graphs.

A graph $G$ is $\alpha ^{+}$-\textit{stable} if $\alpha (G+e)=\alpha (G)$
holds for any edge $e\in E(\overline{G})$, where $\overline{G}$ is the
complement of $G$, \cite{gun}. We shall use the following characterization
that Haynes et al. give for the $\alpha ^{+}$-stable graphs.

\begin{theorem}
\label{th1}\cite{hayn} A graph $G$ is $\alpha ^{+}$-stable if and only if $%
\left| \cap \{S:S\in \Omega (G)\}\right| \leq 1$.
\end{theorem}

Theorem \ref{th1} motivates us to define graph $G$ as $\alpha _{0}^{+}$-%
\textit{stable} if $\left| \cap \{S:S\in \Omega (G)\}\right| =0$, and $%
\alpha _{1}^{+}$-\textit{stable} if $\left| \cap \{S:S\in \Omega
(G)\}\right| =1$, \cite{levm2}. Based on Theorem \ref{th1}, Gunther et al., 
\cite{gun}, give a description of $\alpha ^{+}$\textit{-}stable trees, which
we generalized to bipartite graphs in \cite{levm}. The structure of $\alpha
^{+}$\textit{-}stable bipartite graphs is emphasized in \cite{levm1}.

In this paper we present several properties of K\"{o}nig-Egervary graphs,
which we use further to give necessary and sufficient conditions for
K\"{o}nig-Egervary graphs to be $\alpha ^{+}$\textit{-}stable. We also
characterize K\"{o}nig-Egervary graphs having perfect matchings. Similar
problems related to adding or deleting edges or vertices in connection with
various graph parameters are treated in \cite{acmm}, \cite{brig}, \cite{fav}%
, \cite{ftz}, \cite{mon}, \cite{sumb}.

\section{K\"{o}nig-Egervary Graphs}

Using the definition of K\"{o}nig-Egervary graphs we get:

\begin{lemma}
\label{lem1}($\mathit{i}$) If $G$ is a K\"{o}nig-Egervary graph, then $%
\alpha (G)\geq \left| V(G)\right| /2\geq \mu (G)$.

($\mathit{ii}$) A K\"{o}nig-Egervary graph $G$ has a perfect matching if and
only if $\alpha (G)=\mu (G)$.

($\mathit{iii}$) If $G$ admits a perfect matching, then $\alpha (G)=\mu (G)$
if and only if $G$ is a K\"{o}nig-Egervary graph.
\end{lemma}

For $G_{i},i=1,2,$ let $G=G_{1}*G_{2}$ be the graph with $V(G)=V(G_{1})\cup
V(G_{2})$, and\textit{\ }$E(G)=E(G_{1})\cup E(G_{2})\cup \{xy:$ \textit{for
some} $x\in V(G_{1})\ $\textit{and} $y\in V(G_{2})\}$. Clearly, if $%
H_{1},H_{2}$ are subgraphs of a graph $G$ such that $V(G)=V(H_{1})\cup
V(H_{2})$ and $V(H_{1})\cap V(H_{2})=$ $\emptyset $, then $G=H_{1}*H_{2}$,
i.e., any graph of order at least two admits such decompositions. However,
some particular cases are of special interest. For instance, if: $%
E(H_{i})=\emptyset ,i=1,2,$ then $G=H_{1}*H_{2}$ is bipartite; $%
E(H_{1})=\emptyset $ and $H_{2}$ is complete, then $G=H_{1}*H_{2}$ is a 
\textit{split graph }\cite{hmmr}.

The following proposition shows that the K\"{o}nig-Egervary graphs are, in
this sense, between these two ''extreme'' situations. The equivalence of the
first and the third parts of this result was proposed by Klee without proof
(see \cite{klee}).

\begin{proposition}
\label{prop1}If $G$ is connected, then the following statements are
equivalent:

($\mathit{i}$) $G$ is a K\"{o}nig-Egervary graph;

($\mathit{ii}$) $G=H_{1}*H_{2}$, where $V(H_{1})=S\in \Omega (G)$ and $%
\left| V(H_{1})\right| \geq \mu (G)=\left| V(H_{2})\right| $;

($\mathit{iii}$) $G=H_{1}*H_{2}$, where $V(H_{1})=S$ is a stable set in $G$, 
$\left| S\right| \geq \left| V(H_{2})\right| $, and $(S,V(H_{2}))$ contains
a matching $M$ with $\left| M\right| =\left| V(H_{2})\right| $.
\end{proposition}

\setlength {\parindent}{0.0cm}\textbf{Proof.} ($\mathit{i}$) $\Rightarrow $ (%
$\mathit{ii}$) Let $S\in \Omega (G),H_{1}=G[S]$ and $H_{2}=G[V-S].$ Then we
have $G=H_{1}*H_{2},\alpha (G)+\mu (G)=\left| V(G)\right| =\alpha (G)+\left|
V(H_{2})\right| $, and therefore $\mu (G)=\left| V(H_{2})\right| $. In
addition, Lemma \ref{lem1} ensures that $\left| V(H_{1})\right| \geq \mu (G)$%
.\setlength
{\parindent}{3.45ex}

\ ($\mathit{ii}$) $\Rightarrow $ ($\mathit{iii}$) It is clear if we take the
same $H_{1}$ and $H_{2}$ as in ($\mathit{ii}$).

($\mathit{iii}$) $\Rightarrow $ ($\mathit{i}$) First, we claim that $\left|
M\right| =\mu (G)$. To see this, let assume $W$ be an arbitrary matching in $%
G$ containing some edge of $H$. Since $S$ is stable, we infer that $\left|
W\right| <\left| V(H_{2})\right| =\left| M\right| $. Therefore, $M$ must be
a maximum matching in $G$. Hence we have: $\alpha (G)+\mu (G)\leq \left|
V(G)\right| =\left| S\right| +\left| V(H)\right| =\left| S\right| +\left|
M\right| =\left| S\right| +\mu (G)$, and because $S$ is stable, we obtain
that $\left| S\right| =\alpha (G)$ and $\alpha (G)+\mu (G)=\left|
V(G)\right| $, i.e., $G$ is a K\"{o}nig-Egervary graph. \rule{2mm}{2mm}%
\newline

In the sequel, we shall often represent a K\"{o}nig-Egervary graph $G$ as $%
G=S*H$, where $S\in \Omega (G)$, $H=G[V-S]$, and $\left| V(H)\right| =\mu
(G) $.

\begin{lemma}
\label{lem3}Any maximum matching of a K\"{o}nig-Egervary graph $G=(V,E)$ is
contained in each $(S,V-S)$, where $S\in \Omega (G)$, and, hence, 
\[
\cup \{M:M\ is\ a\ maximum\ matching\ in\ G\}\subseteq \cap \{(S,V-S):S\in
\Omega (G)\}.
\]
\end{lemma}

\setlength {\parindent}{0.0cm}\textbf{Proof.} Let $S\in \Omega (G)$ and $%
G=S*H$. Suppose, on the contrary, that there is a maximum matching $M$ of $G$
and an edge $e=xy\in M\cap E(H).$ Since $S$ is stable, we infer that $\mu
(G)<\left| V(H)\right| $, a contradiction. Therefore, $M$ must be contained
in $(S,V-S)$. \rule{2mm}{2mm}\setlength {\parindent}{3.45ex}\newline

Let $M$ be a maximum matching of a graph $G$. To adopt Edmonds's
terminology, \cite{edm}, we recall the following terms for $G$ relative to $%
M $. The edges in $M$ are \textit{heavy}, while those not in $M$ are \textit{%
light}. An \textit{alternating path} from a vertex $x$ to a vertex $y$ is a $%
x,y$-path whose edges are alternating light and heavy. A vertex $x$ is 
\textit{exposed} relative to $M$ if $x$ is not the endpoint of a heavy edge.
An odd cycle $C$ with $V(C)=\left\{ x_{0},x_{1},...,x_{2k}\right\} \ $and\ $%
E(C)=\left\{ x_{i}x_{i+1}:0\leq i\leq 2k-1\right\} \cup \{x_{2k},x_{0}\}$,
such that $x_{1}x_{2},x_{3}x_{4},...,x_{2k-1}x_{2k}\in M$ is a \textit{%
blossom relative to }$M$\textit{.} The vertex $x_{0}$ is the \textit{base}
of the blossom. The \textit{stem} is an even length alternating path joining
the base of a blossom and an exposed vertex for $M$. The base is the only
common vertex to the blossom and the stem. A \textit{flower} is a blossom
and its stem. A \textit{posy} or a \textit{blossom pair }(cf. \cite{dem})
consists of two (not necessarily disjoint) blossoms joined by an odd length
alternating path whose first and last edges belong to $M$. The endpoints of
the path are exactly the bases of the two blossoms. The following result of
Sterboul, \cite{ster}, characterizes K\"{o}nig-Egervary graphs in terms of
forbidden configurations.

\begin{theorem}
For a graph $G$, the following properties are equivalent:

($\mathit{i}$) $G$ is a K\"{o}nig-Egervary graph;

($\mathit{ii}$) there exist no flower and no posy relative to some maximum
matching $M$;

($\mathit{iii}$) there exist no flower and no posy relative to any maximum
matching $M$.
\end{theorem}

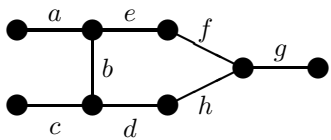
\begin{figure}[h]
\setlength{\unitlength}{1cm}%
\begin{picture}(5,1.6)\thicklines

  \multiput(5,0.3)(1,0){3}{\circle*{0.29}}
  \multiput(5,1.3)(1,0){3}{\circle*{0.29}}
  \multiput(8,0.8)(1,0){2}{\circle*{0.29}}
  \put(5,0.3){\line(1,0){2}}
  \put(5,1.3){\line(1,0){2}}
  \put(8,0.8){\line(1,0){1}}
  \put(6,0.3){\line(0,1){1}}
  \put(7,0.3){\line(2,1){1}}
  \put(7,1.3){\line(2,-1){1}}
\put(5.5,1.5){\makebox(0,0){$a$}}
\put(6.2,0.8){\makebox(0,0){$b$}}
\put(5.5,0){\makebox(0,0){$c$}}

\put(6.5,1.5){\makebox(0,0){$e$}}
\put(6.5,0){\makebox(0,0){$d$}}
\put(7.5,0.3){\makebox(0,0){$h$}}
\put(7.5,1.3){\makebox(0,0){$f$}}
\put(8.5,1){\makebox(0,0){$g$}}
\end{picture}
\caption{Graph $G$ is not blossom-free.}
\label{fig2}
\end{figure}

If a K\"{o}nig-Egervary graph $G$ is blossom-free relative to a maximum
matching $M$, then $G$ is not necessarily blossom-free with respect to any
of its maximum matchings. For instance, the graph $G$ in Figure \ref{fig2}
contains a unique $C_{5}$, which is a blossom relative to the maximum
matching $M_{1}=\left\{ d,e,g\right\} $, and is not a blossom relative to $%
M_{2}=\left\{ a,c,g\right\} $.

\begin{lemma}
\label{lem2}If $M$ is a maximum matching and $S$ is a stable set of a
K\"{o}nig-Egervary graph $G$, then $S\in \Omega (G)$ if and only if $S$
contains all exposed vertices relative to $M$ and one endpoint of each edge
in $M$.
\end{lemma}

\setlength {\parindent}{0.0cm}\textbf{Proof.} According to Proposition \ref
{prop1}, $G=S*H$, where $S\in \Omega (G)$ and $H=G[V-S]$ has $\mu (G)=\left|
V\left( H\right) \right| $. By Lemma \ref{lem3}, $M\subset (S,V-S)$, and
therefore the assertion on $S$ is true.\setlength {\parindent}{3.45ex}

Conversely, since $S$ is stable and $\left| S\right| =\left| M\right|
+\left| V(G)\right| -2\left| M\right| =\left| V(G)\right| -\left| M\right|
=\left| V(G)\right| -\mu (G)$, we get that $S\in \Omega (G)$. \rule{2mm}{2mm}

\begin{theorem}
\label{th4}Let $G$ be a K\"{o}nig-Egervary graph of order at least $2$. Then 
$G$ satisfies $\left| \cap \{V-S:S\in \Omega (G)\}\right| =0$ if and only if
it has a perfect matching and is blossom-free.
\end{theorem}

\setlength {\parindent}{0.0cm}\textbf{Proof.} Suppose that $G$ has no
perfect matching. Let $S\in \Omega (G),G=S*H$, and $M$ be a maximum matching
in $G$. Lemma \ref{lem1} implies that $\alpha (G)>\mu (G)$, and hence, $G$
has at least one exposed vertex $v$ with respect to $M$. Then $v\in S$, and
any $w\in N(v)$ is not contained in $S$. Since the choice of $S$ is
arbitrary, we conclude that $w\in V-S$, for any $S\in \Omega (G)$. Hence, $%
\left| \cap \{V-S:S\in \Omega (G)\}\right| >0$, in contradiction with the
premises on $G$. Thus, $G$ must have a perfect matching. To prove that $G$
is blossom-free, it is sufficient to show that if $G$ is a
K\"{o}nig-Egervary graph, then $\{x:x$\ \textit{is a base of a blossom in}\ $%
G\}$ $\subseteq \cap \{V-S:S\in \Omega (G)\}$, i.e., for any $S\in \Omega (G)
$, no base of a blossom in $G$ belongs to $S$. Let $C$ be a blossom in $G$,
with $V(C)=\left\{ x_{0},x_{1},...,x_{2k}\right\} $, relative to a perfect
matching $M$, and $x_{0}$ be its base. Then $%
x_{1}x_{2},x_{3}x_{4},...,x_{2k-1}x_{2k}\in M$ , and according to Lemma \ref
{lem2}, $S$ contains one of the vertices $x_{1}$ or $x_{2k}$. If $%
x_{1,}x_{2k}\notin S$, then necessarily $x_{2},x_{2k-1}\in S$ and this is
not possible, since the node distance on $C$ between $x_{2}$ and $x_{2k-1}$
is an even number. Hence, $x_{0}\notin S$.\setlength
{\parindent}{3.45ex}

Conversely, Let $M$ be a perfect matching of $G=S*H$, and $b\in V(H)$, where 
$S\in \Omega (G)$, $H=G[V-S]$, and $\left| V(H)\right| =\mu (G)$. We
emphasize a maximum stable set of $G$ that contains $b$. Let denote:

\[
A_{1}=N(b)\cap S,B_{1}=\{b:ab\in M,a\in A_{1}\}, 
\]
\[
A_{2}=N(B_{1})\cap S-A_{1},B_{2}=\{b:ab\in M,a\in A_{2}\}, 
\]
\[
A_{3}=N(B_{2})\cap S-A_{1}\cup A_{2},B_{3}=\{b:ab\in M,a\in A_{3}\},..., 
\]
\[
A_{p}=N(B_{p-1})\cap S-A_{1}\cup ...\cup A_{p-1},B_{p}=\{b:ab\in M,a\in
A_{p}\}, 
\]
and $b\in B=B_{1}\cup ...\cup B_{p},A=A_{1}\cup ...\cup A_{p}$ be such that $%
(B,S-A)=$ $\emptyset .$ Any edge joining two vertices in $B$ would close a
blossom with respect to $M$, which contradicts the fact that $G$ is blossom
free. Therefore, $B$ is stable. The set $B\cup (S-A)$ is also stable,
because $(B,S-A)=$ $\emptyset $. Moreover, $\left| B\right| =\left| A\right| 
$ implies that $B\cup (S-A)$ is a maximum stable set of $G$.

Thus, every $b\in V(H)=V-S$ belongs to a maximum stable set. Since $S$ is
also a maximum stable set, we conclude that any vertex of $G$ belongs to
some maximum stable set of $G$. Clearly, this is equivalent to $\left| \cap
\{V-S:S\in \Omega (G)\}\right| =0$. \rule{2mm}{2mm}\newline

It is worth observing that having a perfect matching is not sufficient for
achieving $\left| \cap \left\{ V-S:S\in \Omega (G)\right\} \right| =0$. For
instance, $K_{4}-e$ is a K\"{o}nig-Egervary graph with perfect matchings,
but $\left| \cap \{V-S:S\in \Omega (G)\}\right| =2$. Being blossom free is
also not enough for $\left| \cap \{V-S:S\in \Omega (G)\}\right| =0$. For
instance, trees without a perfect matching are examples of blossom free
graphs such that $\left| \cap \{V-S:S\in \Omega (G)\}\right| \neq 0$.

\section{$\alpha ^{+}$-Stable K\"{o}nig-Egervary Graphs}

\begin{lemma}
\label{lem4}Any $\alpha ^{+}$-stable K\"{o}nig-Egervary graph has a
near-perfect matching or a perfect matching.
\end{lemma}

\setlength {\parindent}{0.0cm}\textbf{Proof.} Suppose that graph $G$ has
neither a near-perfect matching nor a perfect matching. Let $M$ be a maximum
matching of $G$. Since $\left| V(G)\right| -2\left| M\right| \geq 2$, there
exist two unmatched vertices of $G$, say $x,y$. Hence, $e=xy\in E(\overline{G%
})$, because otherwise $M\cup \{e\}$ is a matching larger than a maximum
matching of $G$. We claim that $x,y$ are contained in all maximum stable
sets of $G$. To see this, let $S\in \Omega (G)$ and $H=G[V-S].$ Then $G=S*H$%
, where $\mu (G)=\left| V(H)\right| =\left| M\right| $. By Lemma \ref{lem3},
we have that $M\subseteq (S,V-S)$. Hence, $x,y$ $\in S$, because these
vertices are unmatched and non-adjacent. Since $S$ was an arbitrary maximum
stable set of $G$, we infer that $x,y\in \cap \{S:S$ $\in \Omega (G)\}$. By
Theorem \ref{th1}, it contradicts the fact that $G$ is $\alpha ^{+}$-stable.
Consequently, $G$ must have a near-perfect matching or a perfect matching. 
\rule{2mm}{2mm}\setlength {\parindent}{3.45ex}

\begin{theorem}
\label{th22}A K\"{o}nig-Egervary graph $G$ is $\alpha ^{+}-$stable if and
only if it has a perfect matching and $\left| \cap \{V-S:S\in \Omega
(G)\}\right| \leq 1$.
\end{theorem}

\setlength {\parindent}{0.0cm}\textbf{Proof.} Let $G$ be $\alpha ^{+}-$%
stable, and $S$ $\in \Omega (G)$. Suppose, on the contrary, that $G$ has no
perfect matching, i.e., by Lemma \ref{lem1}, $\alpha (G)>\mu (G)$. Lemma \ref
{lem4} implies that $G$ has a near-perfect matching $M$, which is contained,
according to Lemma \ref{lem3}, in $(S,V-S).$ Hence, we get that $\alpha
(G)=\left| S\right| =\mu (G)+1=\left| V-S\right| +1=\left| M\right| +1$, and
there are $x,y\in S$ and $z\in V-S$ such that $xz\in E(G)-M$ and $yz\in M$.
We claim that $x,y$ belong also to any other maximum stable set $W$ of $G$,
since otherwise if:\setlength {\parindent}{3.45ex}

\textit{(}$\mathit{a}$\textit{)} $z\in W,$ then $x,y\notin W$, and hence $%
\left| W\right| <\alpha (G),$ a contradiction;

\textit{(}$\mathit{b}$\textit{)} only $x\in W$ or only $y\in W$, then $%
z\notin W$, and again the contradiction $\left| W\right| <\alpha (G)$,
because all vertices of $S-\{x,y\}$ are respectively matched, by $M\emph{-}%
\{yz\}$, with vertices in $V-S-\{z\}$.

Thus, we get that $x,y\in \cap \{S:S$ $\in \Omega (G)\}$, and according to
Theorem \ref{th1}, this contradicts the fact that $G$ is $\alpha ^{+}$%
-stable. Therefore, $G$ has a perfect matching, say $M$. Then, for any edge $%
e=xy\in M$, we have that $x\in \cap \{S:S\in \Omega (G)\}$\ if and only if $%
y\in \cap \{V-S:S\in \Omega (G)\}$. Consequently, we obtain that $\left|
\cap \{S:S\in \Omega (G)\}\right| =\left| \cap \{V-S:S\in \Omega
(G)\}\right| $, and Theorem \ref{th1} implies $\left| \cap \{V-S:S\in \Omega
(G)\}\right| \leq 1$.

Conversely, suppose $G$ has a perfect matching and $\left| \cap \{V-S:S\in
\Omega (G)\}\right| \leq 1$. As we saw in the previous paragraph, the
existence of a perfect matching in $G$ results in $\left| \cap \{S:S\in
\Omega (G)\}\right| =\left| \cap \{V-S:S\in \Omega (G)\}\right| $. Since $%
\left| \cap \{V-S:S\in \Omega (G)\}\right| \leq 1$, Theorem \ref{th1}
ensures that $G$ is $\alpha ^{+}$-stable. \rule{2mm}{2mm}\newline

It is worth mentioning that there are K\"{o}nig-Egervary graphs with perfect
matchings, which are not $\alpha ^{+}$-stable; e.g., the graph $K_{4}-e$.
However, for bipartite graphs, this condition is also sufficient (see
Corollary \ref{cor1}).

Propositions \ref{prop5} and \ref{prop6} show that any $\alpha ^{+}$-stable
K\"{o}nig-Egervary graph $G$ with $\left| \cap \{V-S:S\in \Omega
(G)\}\right| =1$ may be decomposed into two $\alpha ^{+}$-stable
K\"{o}nig-Egervary graphs $G_{1},G_{2}$ with $\left| \cap \{V-S:S\in \Omega
(G_{1})\}\right| =0$, and $\left| \cap \{V-S:S\in \Omega (G_{2})\}\right| =0$%
.

\begin{proposition}
\label{prop5}If $G$ is a K\"{o}nig-Egervary graph with $\left| \cap
\{V-S:S\in \Omega (G)\}\right| =1$ and $\alpha (G)=\mu (G)$, then there
exists $xy\in E(G),$ such that $H=G-\{x,y\}$ is a K\"{o}nig-Egervary graph
with\ $\left| \cap \{V-S:S\in \Omega (H)\}\right| =0$ and $\alpha (H)=\mu (H)
$.
\end{proposition}

\setlength {\parindent}{0.0cm}\textbf{Proof.} Let $M$ be a perfect matching
of $G$, which exists by Lemma \ref{lem1}. Let also $x\in \cap \{V-S:S\in
\Omega (G)\}$ and $y\in V(G)$ be such that $e=xy\in M$. Hence, it follows
that $y\in \cap \{S:S\in \Omega (G)\}$, and therefore, $H=G-\{x,y\}$ is a
K\"{o}nig-Egervary graph with $\alpha (H)=\mu (H)\ $and$\mathit{\ }\left|
\cap \{V-S:S\in \Omega (H)\}\right| =0$. \rule{2mm}{2mm}%
\setlength
{\parindent}{3.45ex}

\begin{proposition}
\label{prop6}If $G$ is a K\"{o}nig-Egervary graph with $\left| \cap
\{V-S:S\in \Omega (G)\}\right| =0$ and $K_{2}=\{\{x,y\},\{xy\}\}$, then
every graph $F=G+K_{2}$ having: 
\[
V(F)=V(G)\cup \{x,y\},E(F)\supseteq E(G)\cup \{xy\}\ and\ \{(y,S)\neq
\emptyset ,\ for\ all\ S\in \Omega (G)\},
\]
is a K\"{o}nig-Egervary graph with a perfect matching, and\textit{\ }$\left|
\cap \{V-S:S\in \Omega (F)\}\right| =1$.
\end{proposition}

\setlength {\parindent}{0.0cm}\textbf{Proof.} By Theorem \ref{th4}, we
obtain that $G$ admits perfect matchings. Since $(y,S)\neq \emptyset $ for
any $S\in \Omega (G)$, we get $\Omega (F)=\{S\cup \{x\}:S\in \Omega (G)\}$. $%
M\cup \{xy\}$ is a perfect matching in $F$, for any perfect matching $M$ of $%
G$. Consequently, $F$ has a perfect matching and $\cap \left\{ V-S:S\in
\Omega (F)\right\} =\{y\}$. According to Proposition \ref{prop1}, $F$ is
also a K\"{o}nig-Egervary graph. \rule{2mm}{2mm}%
\setlength
{\parindent}{3.45ex}\newline

The next theorem presents a more specific characterization of $\alpha ^{+}$%
-stable K\"{o}nig-Egervary graphs.

\begin{theorem}
\label{th8}If $G$ is a K\"{o}nig-Egervary graph of order at least $2$, then
the following statements are equivalent:

($\mathit{i}$) $G$ is $\alpha ^{+}$-stable;

($\mathit{ii}$) either $\left| \cap \{V-S:S\in \Omega (G)\}\right| =0$, or $%
\left| \cap \{V-S:S\in \Omega (G)\}\right| =1$, and $G$ has a perfect
matching;

($\mathit{iii}$) $G$ has a perfect matching, and either there exists $xy\in
E(G)$, such that $H=G-\{x,y\}$ is blossom-free and has a perfect matching,
or $G$ is blossom-free.
\end{theorem}

\setlength {\parindent}{0.0cm}\textbf{Proof.} ($\mathit{iii}$) If $\left|
\cap \{V-S:S\in \Omega (G)\}\right| =0$, then $G$ has a perfect matching and
it is blossom-free, by Theorem \ref{th4}. If $\left| \cap \{V-S:S\in \Omega
(G)\}\right| =1$ and $G$ has a perfect matching, then Proposition \ref{prop5}
and Lemma \ref{lem1}, imply that there exists $xy\in E(G)$, such that $%
H=G-\{x,y\}$ is blossom-free and has a perfect matching.%
\setlength
{\parindent}{3.45ex}

($\mathit{iii}$) $\Rightarrow $ ($\mathit{i}$) If $G$ has a perfect matching
and is blossom-free, then according to Theorem \ref{th4}, we get $\left|
\cap \{V-S:S\in \Omega (G)\}\right| =0$, and further, Theorem \ref{th22}
ensures that $G$ is $\alpha ^{+}$-stable. If $G$ has a perfect matching, and
there exist $xy\in E(G)$, such that $H=G-\{x,y\}$ is blossom-free and has a
perfect matching, then $G$ is $\alpha ^{+}$-stable, according to Proposition 
\ref{prop6}. \rule{2mm}{2mm}\newline

The graph $K_{4}-e$ shows that it is not enough to have a perfect matching
in order to ensure that a K\"{o}nig-Egervary graph is $\alpha ^{+}$-stable.

Notice also that $P_{3}$ is a K\"{o}nig-Egervary graph, $\left| \cap
\{V-S:S\in \Omega (P_{3})\}\right| =1$, but $P_{3}$ is not $\alpha ^{+}$%
-stable.

\begin{figure}[h]
\setlength{\unitlength}{1cm}%
\begin{picture}(5,2)\thicklines
  \multiput(5,0.5)(1,0){4}{\circle*{0.29}}
  \multiput(5,1.5)(1,0){4}{\circle*{0.29}}
  \put(5,0.5){\line(1,0){3}}
  \put(5,1.5){\line(1,0){3}}
  \put(5,0.5){\line(1,1){1}}
  \put(6,0.5){\line(1,1){1}}
  \put(7,0.5){\line(1,1){1}}
  \put(5,1.5){\line(1,-1){1}}
  \put(6,1.5){\line(1,-1){1}}
  \put(7,1.5){\line(1,-1){1}}
  \put(6,0.5){\line(0,1){1}}
  \put(7,0.5){\line(0,1){1}}
\put(4.7,0.5){\makebox(0,0){$v_{5}$}}
\put(4.7,1.5){\makebox(0,0){$v_{1}$}}
\put(6,0){\makebox(0,0){$v_{6}$}}
\put(6,1.9){\makebox(0,0){$v_{2}$}}
\put(7,0){\makebox(0,0){$v_{7}$}}
\put(7,1.9){\makebox(0,0){$v_{3}$}}
\put(8.4,0.5){\makebox(0,0){$v_{8}$}}
\put(8.4,1.5){\makebox(0,0){$v_{4}$}}
\put(5.5,0.2){\makebox(0,0){$c$}}
\put(5.5,1.7){\makebox(0,0){$a$}}
\put(7.5,0.2){\makebox(0,0){$d$}}
\put(7.5,1.7){\makebox(0,0){$b$}}
\end{picture}
\caption{$G$ is not-blossom-free, $H=G-\{v_{2},v_{3}\}$ is blossom-free, and 
$\alpha (H)\neq \mu (H)$.}
\label{fig3}
\end{figure}
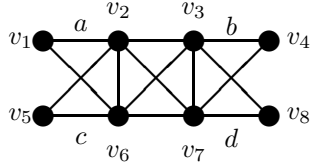

Observe that the graph $G$, in Figure \ref{fig3}, has blossoms with respect
to the perfect matching $M=\left\{ a,b,c,d\right\} $, and only for $x\in
\left\{ v_{2},v_{6}\right\} $ and $y\in \left\{ v_{3},v_{7}\right\} $, the
corresponding subgraph $H=G-\{x,y\}$ is connected and blossom-free, but $%
\alpha (H)\neq \mu (H)$. In addition, $G$ is a K\"{o}nig-Egervary non-$%
\alpha ^{+}$-stable graph, since $\alpha (G+v_{1}v_{5})=3<\alpha (G)$.

\section{Applications of $\alpha ^{+}$-stable K\"{o}nig-Egervary Graphs}

Combining Theorems \ref{th8}, \ref{th1} and Lemma \ref{lem1}, we obtain:

\begin{corollary}
\label{cor1}\cite{levm} If $G$ is a bipartite graph, then the following
assertions are equivalent:

\textit{(}$\mathit{i}$\textit{)} $G$ is $\alpha ^{+}$-stable;

\textit{(}$\mathit{ii}$\textit{) }$G$ possesses a perfect matching;

\textit{(}$\mathit{iii}$\textit{)} $G$ has two maximum stable sets that
partition its vertex set;

\textit{(}$\mathit{iv}$\textit{)} $\left| \cap \{S:S\in \Omega (G)\}\right|
=0$.
\end{corollary}

\begin{figure}[h]
\setlength{\unitlength}{1cm}%
\begin{picture}(5,1.5)\thicklines

  \multiput(2,0)(1,0){3}{\circle*{0.29}}
  \multiput(2,1)(1,0){3}{\circle*{0.29}}
  \multiput(5,0.5)(1,0){2}{\circle*{0.29}}
  \put(5,0.5){\line(1,0){1}}
  \put(3,0){\line(0,1){1}}
  \put(4,0){\line(0,1){1}}
  \put(4,0){\line(2,1){1}}
  \put(4,1){\line(2,-1){1}}
  \put(2,0){\line(1,1){1}}
  \put(2,1){\line(1,-1){1}}
  \put(2,0){\line(2,1){2}}
  \put(2,1){\line(2,-1){2}}
  \put(3,0){\line(1,1){1}}
  \put(3,1){\line(1,-1){1}}
\put(1.5,0.5){\makebox(0,0){$G_{1}$}}

 \multiput(7,0)(1,0){4}{\circle*{0.29}}
 \multiput(8,1)(1,0){2}{\circle*{0.29}}
  \put(7,0){\line(1,0){3}}
  \put(8,1){\line(1,0){1}}
  \put(8,0){\line(1,1){1}}
  \put(9,0){\line(0,1){1}}
  \put(10.5,0.5){\makebox(0,0){$G_{2}$}}
\end{picture}
\caption{$\alpha ^{+}$-stable non-bipartite K\"{o}nig-Egervary graphs.}
\label{fig4}
\end{figure}
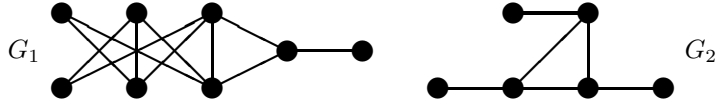

In other words, the bipartite graphs can be only $\alpha _{0}^{+}$-stable.
Nevertheless, there exist non-bipartite K\"{o}nig-Egervary $\alpha _{0}^{+}$%
-stable graphs (e.g., $G_{2}$ in Figure \ref{fig4}), and also non-bipartite
K\"{o}nig-Egervary $\alpha _{1}^{+}$-stable graphs (e.g., $G_{1}$ in Figure 
\ref{fig4}).

\begin{proposition}
\label{construction}Let $G$ be an $\alpha ^{+}$-stable bipartite graph. If
graph $H=G\bullet K_{p}$ has $V(H)=V(G)\cup V(K_{p})$ and $E(H)=E(G)\cup
E(K_{p})\cup W$, where:

\textit{(}$\mathit{i}$\textit{)} $W=\{xa,xb\}$, with $ab$ in a perfect
matching of $G$, for $p\leq 2$ and $x\in V(K_{p})$; or

\textit{(}$\mathit{ii}$\textit{) }$W=\{xy\}$, for some $x\in V(K_{p})$ and $%
y\in V(G),$ for $p\geq 3$,

then $H$ is $\alpha ^{+}$-stable.
\end{proposition}

\setlength {\parindent}{0.0cm}\textbf{Proof.} \textit{(}$\mathit{i}$\textit{)%
} If $p=1$, we claim that $\alpha (H)=\alpha (G)$. Otherwise, there is a
stable set $S$ in $H$ with $\left| S\right| =\alpha (H)>\alpha (G)$, and
consequently $S-\left\{ x\right\} $ is a maximum stable set in $G$ that
contains neither $a$, nor $b$, a contradiction, since $G$ is $\alpha ^{+}$%
-stable. Hence, we have $\Omega (H)=\Omega (G)$ and clearly, $\cap \{S:S\in
\Omega (H)\}=\emptyset $. If $p=2$, then $\left| \cap \{S:S\in \Omega
(H)\}\right| =\left| \{y\}\right| 1$, because $G$ is $\alpha ^{+}$-stable
and any maximum stable set of $H$ is of the form $S\cup \{y\}$, where $S\in
\Omega (G)$ and $y\in V(K_{2})-\{x\}$. Therefore, by Theorem \ref{th1}, $H$
is $\alpha ^{+}$-stable.\setlength
{\parindent}{3.45ex}

\textit{(}$\mathit{ii}$\textit{)} In this case, $\cap \{S:S\in \Omega
(H)\}=\emptyset $, because $G$ is $\alpha ^{+}$-stable and any maximum
stable set of $H$ is of the form $S\cup \{z\}$, where $S\in \Omega (G)$ and $%
z\in V(K_{p})-\{x\}$. According to Theorem \ref{th1}, $H$ is $\alpha ^{+}$%
-stable. \rule{2mm}{2mm}\newline

The graph $G=K_{1}\bullet K_{p}$, for $p\geq 4$, $K_{1}=(\left\{ x\right\}
,\emptyset ),V(G)=\left\{ x\right\} \cup V(K_{p})$ and $E(G)=\{xy\}\cup
E(K_{p})$, where $y\in V(K_{p})$, is $\alpha ^{+}$-stable,
non-K\"{o}nig-Egervary graph, and $\left| \cap \{S:S\in \Omega (G)\}\right|
=\left| \left\{ x\right\} \right| =1$. Taking also into account Proposition 
\ref{construction}, we obtain the following:

\begin{corollary}
For every natural number $n\geq 5$ there exist $\alpha ^{+}$-stable
non-K\"{o}nig-Egervary graphs $G_{1},G_{2}$ of order $n$ such that 
\[
\left| \cap \{S:S\in \Omega (G_{1})\}\right| =0\ and\ \left| \cap \{S:S\in
\Omega (G_{2})\}\right| =1.
\]
\end{corollary}

\begin{proposition}
\label{semi}If $G$ is a K\"{o}nig-Egervary graph of order $n\geq 2$, and $%
\alpha (G)>n/2$, then $\left| \cap \{S:S\in \Omega (G)\}\right| \geq 2$.
\end{proposition}

\setlength {\parindent}{0.0cm}\textbf{Proof.} If $\alpha (G)>n/2$ , then $G$
has no perfect matching, and by Theorem \ref{th22}, $G$ is not $\alpha ^{+}$%
-stable. Consequently, Theorem \ref{th1} implies that $\left| \cap \{S:S\in
\Omega (G)\}\right| \geq 2$. \rule{2mm}{2mm}\setlength
{\parindent}{3.45ex}\newline

For general case, it has been proven that:

\begin{proposition}
\label{semi1}\cite{hamhansim} If $G$ is a graph with $\alpha (G)>\left|
V(G)\right| /2$ , then 
\[
\left| \cap \{S:S\in \Omega (G)\}\right| \geq 1.
\]
For example, in a bipartite graph $G=(A,B,E)$ such that $\left| A\right|
\neq \left| B\right| $, there exists at least one vertex belonging to all
maximum stable sets of $G$, i.e., 
\[
\left| \cap \{S:S\in \Omega (G)\}\right| \geq 1.
\]
\end{proposition}

Since any bipartite graph is K\"{o}nig-Egervary, Proposition \ref{semi}
yields the following result, which has been already done independently in 
\cite{levm1}, as a strengthening of Proposition \ref{semi1} in the case of
bipartite graphs.

\begin{corollary}
\label{cor6}If $G=(A,B,E)$ is bipartite and $\left| A\right| \neq \left|
B\right| $, then 
\[
\left| \cap \{S:S\in \Omega (G)\}\right| \geq 2.
\]
\end{corollary}

A \textit{pendant edge} is an edge incident with a \textit{pendant vertex}
(i.e., a vertex of degree one). A vertex $v$ is $\alpha $\textit{-critical}
in $G$ if $\alpha (G-v)<\alpha (G)$.

\begin{theorem}
\label{th9}For a graph $G$ of order at least two, the following are
equivalent:

($\mathit{i}$) $G$ has a perfect matching $M$ consisting of its pendant
edges;

($\mathit{ii}$) $G$ has exactly $\alpha (G)$ pendant vertices and none of
them is $\alpha $-critical;

($\mathit{iii}$) $G$ is a K\"{o}nig-Egervary $\alpha ^{+}$-stable graph with
exactly $\alpha (G)$ pendant vertices.
\end{theorem}

\setlength {\parindent}{0.0cm}\textbf{Proof.} ($\mathit{i}$) $\Rightarrow $ (%
$\mathit{ii}$) It is clear that $S=\{x:x$\ \textit{is a pendant vertex in}\ $%
G\}$ is stable in $G$. If $\alpha (G)>\left| S\right| =\left| V(G)\right| /2$%
, then any maximum stable set $W$ of $G$ must contain some pair of vertices,
matched by $M$, a contradiction, since $W$ is stable. Hence, $\left|
S\right| =\alpha (G)$ holds. In addition, if $x\in S$ and $y$ is its single
neighbor in $G$, then $S\cup \{y\}-\left\{ x\right\} $ is a maximum stable
set in $G-\left\{ x\right\} $, i.e., $x$ is not $\alpha $-critical in $G$.%
\setlength {\parindent}{3.45ex}

($\mathit{ii}$) $\Rightarrow $ ($\mathit{iii}$) Now, clearly $S=\{x:x$\ 
\textit{is a pendant vertex in}\ $G\}\in \Omega (G)$ and let denote $%
M=\{xy:x\in S$\ \textit{and} $y\in N(x)\}$. By Proposition \ref{prop1}, $%
G=S*H$ and if some $z\in V(H)$ is not matched by $M$, then $S\cup \left\{
z\right\} $ is a stable set larger than $S$, a contradiction. Hence, we get
that $\left| M\right| =\mu (G)=\left| V(H)\right| $, i.e., $G$ is a
K\"{o}nig-Egervary graph. According to Theorem \ref{th8}, $G$ is also $%
\alpha ^{+}$-stable, because is blossom-free with respect to $M$ and another
perfect matching does not exist.

($\mathit{iii}$) $\Rightarrow $ ($\mathit{i}$) According to Theorem \ref{th8}%
, $G$ has a perfect matching $M$, and since $\{x:x$\ \textit{is a pendant
vertex in}\ $G\}\in \Omega (G)$, $M$ consists of all the pendant edges of $G$%
. \rule{2mm}{2mm}

\begin{lemma}
\label{lem5}If $G$\ is a K\"{o}nig-Egervary graph, then 
\[
N(\cap \{S:S\in \Omega (G)\})=\cap \{V-S:S\in \Omega (G)\}.
\]
\end{lemma}

\setlength {\parindent}{0.0cm}\textbf{Proof.} Let denote $A=\cap \{S:S\in
\Omega (G)\}$ and $B=\cap \{V-S:S\in \Omega (G)\}$. If $v\in N(A)$, then
clearly $v\notin S$, for any $S\in \Omega (G)$, i.e., $N(A)\subseteq B$. Let 
$M$ be a maximum matching of $G$ and $x\in B$. According to Lemma \ref{lem3}%
, $M\subset (S,V(G)-S)$ holds for any $S\in \Omega (G)$, and by Proposition 
\ref{prop1}, we have also $\left| M\right| =\left| V(G)-S\right| $. Since $%
x\in B$, it follows that there is $xy\in M$, and hence, Lemma \ref{lem2}
implies that $y\in S$, for any $S\in \Omega (G)$, i.e., $y\in A$.
Consequently, we get that $x\in N(A)$, and because $x$ was an arbitrary
vertex of $B$, it results $B\subseteq N(A)$, and this completes the proof. 
\rule{2mm}{2mm}\setlength {\parindent}{3.45ex}\newline

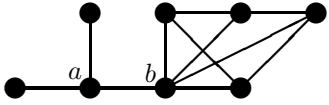
\begin{figure}[h]
\setlength{\unitlength}{1cm}%
\begin{picture}(5,1.5)\thicklines

  \multiput(4,0)(1,0){4}{\circle*{0.29}}
  \multiput(5,1)(1,0){4}{\circle*{0.29}}
  \put(4,0){\line(1,0){3}}
  \put(6,1){\line(1,0){2}}
  \put(5,0){\line(0,1){1}}
  \put(6,0){\line(0,1){1}}
  \put(6,0){\line(1,1){1}}
  \put(6,1){\line(1,-1){1}}
  \put(6,0){\line(2,1){2}}
  \put(7,0){\line(1,1){1}}
  \put(4.8,0.2){\makebox(0,0){$a$}}
  \put(5.8,0.2){\makebox(0,0){$b$}}
\end{picture}
\caption{$N(\cap \{S:S\in \Omega (G)\})=\{a\}\neq \{a,b\}=\cap \{V-S:S\in
\Omega (G)$, and $G$ is a non-K\"{o}nig-Egervary graph without perfect
matchings.}
\label{fig5}
\end{figure}

Lemma \ref{lem5} is not true for general graphs; e.g., the graph in Figure 
\ref{fig5}.

\begin{lemma}
\label{lem6}If $G$\ is a K\"{o}nig-Egervary graph and $M$ is a maximum
matching, then $M$ matches $N(\cap \{S:S\in \Omega (G)\})$ into $\cap
\{S:S\in \Omega (G)\}$.
\end{lemma}

\setlength {\parindent}{0.0cm}\textbf{Proof.} In accordance with Proposition 
\ref{prop1}, $G$ can be written as $G=S*H$, where $S\in \Omega (G)$, $%
H=(V(H),E(H))=G[V-S]$, and $\left| V(H)\right| =\mu (G)$. By Lemma \ref{lem3}%
, $M\subset (S,V(G)-S)=(S,V(H))$, and, clearly, $N(\cap \{S:S\in \Omega
(G))\subset V(G)-S$. Hence, any $x\in N(\cap \{S:S\in \Omega (G))$ is
matched with some $y\in S$. Moreover, according to Lemma \ref{lem2}, if $x$
belongs to no maximum stable set of $G$, then $y\in \cap \{S:S\in \Omega
(G)\}$. Therefore, $M$ matches $N(\cap \{S:S\in \Omega (G))$ into $\cap
\{S:S\in \Omega (G)$. \rule{2mm}{2mm}\setlength {\parindent}{3.45ex}

\begin{theorem}
\label{th10}If $G$\ is a K\"{o}nig-Egervary graph, then $G$ has a perfect
matching if and only if $\left| \cap \{S:S\in \Omega (G)\}\right| =\left|
\cap \{V-S:S\in \Omega (G)\}\right| $.
\end{theorem}

\setlength {\parindent}{0.0cm}\textbf{Proof.} Let $M$ be a perfect matching
of $G$. Then, for any edge $e=xy\in M$, we have that $x\in \cap \{S:S\in
\Omega (G)\}$\ if and only if $y\in \cap \{V-S:S\in \Omega (G)\}$.
Consequently, we get that $\left| \cap \{S:S\in \Omega (G)\}\right| =\left|
\cap \{V-S:S\in \Omega (G)\}\right| $.\setlength
{\parindent}{3.45ex}

Conversely, assume that $\left| \cap \{S:S\in \Omega (G)\}\right| =\left|
\cap \{V-S:S\in \Omega (G)\}\right| $. Let $A=\cap \{S:S\in \Omega
(G)\},B=\cap \{V-S:S\in \Omega (G)\}$, and $G_{0}=G-N[A]$. By Proposition 
\ref{prop1}, $G=S*H$, where $H=(V(H),E(H))=G[V-S]$ has $\left| V(H)\right|
=\mu (G)$, and $S\in \Omega (G)$. If $M$ is a maximum matching, then by
Lemma \ref{lem6}, $M$ matches $N(A)$ into $A$. Hence, $\left| A\right| \geq
\left| N(A)\right| $. Since, by Lemma \ref{lem5}, $N(A)=B$, we get that $%
\left| A\right| \geq \left| N(A)\right| =\left| B\right| $, and consequently 
$\left| A\right| =\left| N(A)\right| $. Therefore the restriction $M_{1}$ of 
$M$ on $G[A\cup N\left( A\right) ]$ is a perfect matching.

For a K\"{o}nig-Egervary graph $G,\alpha \left( G_{0}\right) =\alpha \left(
G\right) -\left| A\right| ,$ $\mu \left( G_{0}\right) =\mu \left( G\right)
-\left| N(A)\right| $. In our case, when $\left| A\right| =\left|
N(A)\right| $, it means that $G_{0}$ is a K\"{o}nig-Egervary graph, as well.
Moreover, $\left| \cap \{S:S\in \Omega (G_{0})\}\right| =0$, and
consequently according to Theorem \ref{th1} $G_{0}$ is an $\alpha ^{+}$%
-stable graph. By Theorem \ref{th22}, $G_{0}$ has a perfect matching, say $%
M_{0}$, which together with $M_{1}$ builds a perfect matching of $G$. \rule%
{2mm}{2mm}\newline

It is interesting to mention that there exist non-K\"{o}nig-Egervary graphs
enjoying the equality $\left| \cap \{S:S\in \Omega (G)\}\right| =\left| \cap
\{V-S:S\in \Omega (G)\}\right| $ without perfect matchings (e.g., the graph
in Figure \ref{fig5}).

Combining Theorem \ref{th10} and Corollary \ref{cor1} we obtain:

\begin{corollary}
If $G$ is bipartite, then $\left| \cap \{S:S\in \Omega (G)\}\right| =\left|
\cap \{V-S:S\in \Omega (G)\}\right| $ if and only if $\left| \cap \{S:S\in
\Omega (G)\}\right| =\left| \cap \{V-S:S\in \Omega (G)\}\right| =0.$
\end{corollary}

\section{Conclusions and future work}

In this paper we return the attention of the reader to the notion of a
K\"{o}nig-Egervary graph. We state several properties of K\"{o}nig-Egervary
graphs, showing that these graphs give a fruitful developing of the
bipartite graphs theory. Our main findings refer to the $\alpha ^{+}$%
-stability of K\"{o}nig-Egervary graphs. These results generalize some
previously known statements for trees and bipartite graphs. In addition, we
characterize those K\"{o}nig-Egervary graphs for which ''to be blossom-free
relative to some perfect matching'' is equivalent to ''to be blossom-free
relative to any perfect matching''. This condition is similar both in form
and spirit to Sterboul's characterization of K\"{o}nig-Egervary graphs. An
obvious question arises: which K\"{o}nig-Egervary graphs are $\alpha ^{-}$%
-stable (i.e., have stability number insensitive to deletion of any edge)?
It would be also interesting to describe the K\"{o}nig-Egervary graphs that
are both $\alpha ^{-}$-stable and $\alpha ^{+}$-stable.

\end{document}